\documentclass[reqno]{amsart}
\usepackage{amsthm,amsmath,amssymb,amsfonts}
\usepackage{epsfig}
\usepackage{graphicx,color}
\usepackage{pstricks}
\usepackage{verbatim}
\usepackage[T1]{fontenc}
\usepackage[ansinew]{inputenc}
\usepackage{mathrsfs}
\usepackage[backref,colorlinks]{hyperref}

\newtheorem{theorem}{Theorem}

\newtheorem{lemma}[theorem]{Lemma}

\newtheorem{claim}[theorem]{Claim}

\newtheorem{obs}[theorem]{Observation}
\newtheorem{conjecture}[theorem]{Conjecture}

\let\theta=\vartheta
\let\rho=\varrho
\let\sigma=\varsigma

\let\polishlcross=\l
\def\l{\ifmmode\ell\else\polishlcross\fi}

\def\cF{{\mathcal F}}

\def\cP{{\mathcal P}}

\def\l{\lambda}

\newcommand{\proofstart}{\noindent {\bf Proof:\hspace{2em}}}

\newcommand{\proofend}{\hspace*{\fill}\mbox{$\Box$}\vspace*{4mm}}

\newcommand{\rt}{\right}

\newcommand{\lt}{\left}

\begin{document}
\title[How many colors guarantee a rainbow matching?]{How many colors guarantee a rainbow matching?}

\author[R.~Glebov]{Roman Glebov}
\address{Institut f\"{u}r Mathematik, Freie Universit\"at Berlin, Arnimallee 3-5, D-14195 Berlin, Germany}
\email{glebov@mi.fu-berlin.de}
\thanks{The first author was supported by DFG within the research training group "Methods for Discrete Structures".}
\author[B.~Sudakov]{Benny Sudakov}
\address{Department of Mathematics, UCLA, Los Angeles, CA 90095}
\email{bsudakov@math.ucla.edu}
\thanks{The second author was supported
in part by NSF grant DMS-1101185, by AFOSR MURI grant FA9550-10-1-0569 and by a USA-Israel BSF
grant.}
\author[T.~Szab\'o]{Tibor Szab\'o}
\address{Institut f\"{u}r Mathematik, Freie Universit\"at Berlin, Arnimallee 3-5, D-14195 Berlin, Germany}
\email{szabo@mi.fu-berlin.de}

\date{}

\begin{abstract}
Given a coloring of the  edges of a multi-hypergraph, a rainbow $t$-matching is a collection of 
$t$ disjoint edges, each having a different color.
In this note we study the problem 
of finding a rainbow $t$-matching in an $r$-partite $r$-uniform multi-hypergraph whose edges are colored 
with $f$ colors such that every color class is a matching of size $t$.
This problem was posed by Aharoni and Berger~\cite{AB09}, who asked to determine the minimum number of 
colors which guarantees a rainbow matching. We improve on the known upper bounds for this problem for all values of the parameters.
In particular for every fixed $r$, we give an upper bound which is polynomial in $t$,
improving the superexponential estimate of Alon~\cite{noga}.
Our proof also works in the setting not requiring the hypergraph to be  $r$-partite.
\end{abstract}

\maketitle

\section{Introduction}

An $r$-{\em uniform multi-hypergraph} is a pair $(V,E)$, where $V$ is the set of vertices of $G$ and
$E$ is a multiset of
$r$-element subsets of the vertex set $V$.
In this paper we usually refer to $r$-uniform multi-hypergraphs as {\em $r$-graph}s.
A {\em matching} is a set of pairwise disjoint edges.
Given a coloring $c: E(G)\rightarrow [f]$ of the edges of an $r$-graph $G$, we call a matching $M\subseteq E(G)$
a {\em rainbow matching} if all its edges have distinct colors.
An $r$-graph $G=(V,E)$ is called {\em $r$-partite} if the vertex set $V=V_1\cup \cdots \cup V_r$
is the disjoint union of $r$ parts $V_i$ and every edge $e\in E$ intersects each part in exactly one
vertex, i.e., $|e\cap V_i|=1$ for all $e\in E$ and $i=1, \ldots , r$.

The problem of finding rainbow subgraphs in edge-colored graphs/hypergraphs has long history and goes back more than 60 years 
to Canonical Ramsey Theorem of Erd\H{o}s and Rado \cite{ER50}. One particular setting, which was extensively studied and which also appears naturally
in problems of additive combinatorics, is when the edge-coloring of the host graph/hypergraph is {\em proper}, i.e., 
every color class in the coloring forms a matching  (see e.g. \cite{Ba, AJMP, 
KMSV} and their references). In this paper we also consider properly edge-colored hypergraphs in which we want to find
a large rainbow matchings. Being interesting for its own right, the problem of finding rainbow matchings in hypergraphs 
can also be used to study various classical extremal problems. For example, consider the following old conjecture of Ryser.
A subset $U\subseteq V(G)$ of the vertices forms a {\em vertex cover} of multi-hypergraph $G$ if every edge of $G$ intersects $U$.
The {\em covering number $\tau(G)$} is the size of the smallest vertex cover of $G$ and   
the {\em matching number} $\nu (G)$ is the size of the largest matching in $G$.
Since the union of the edges of a maximum matching in $G$ is a vertex cover, for every $r$-graph $G$ we have that
\begin{align} \label{eq:tau1} \nu (G) & \geq \tau (G)/r. \end{align}
Ryser's Conjecture states that this trivial lower bound on $\nu (G)$ can be improved to $\tau(G)/(r-1)$ provided $G$ is $r$-partite.
The conjecture follows for $r=2$ from K\H onig's theorem, was proved by Aharoni~\cite{A01} for $r=3$, and is still open for all $r\geq 4$. 
One approach to this problem is to 
specify one of the vertex classes, say $V_1$, and consider the union $H$ of the links of each vertex $x\in V_1$.
(The link of a vertex $x\in V_1$ is the $(r-1)$-graph consisting of those $(r-1)$-sets in
$V_2 \cup \cdots \cup V_r$ which, together with $x$,  form an edge of $G$.)
Then $H$ is an $(r-1)$-partite $(r-1)$-graph, whose edges are colored by
the elements of $V_1$: an edge of $H$ coming from the link of some vertex $x\in V_1$ is colored with ''color'' $x$.
Looking for  a matching in the original $r$-graph $G$ corresponds exactly to
looking for a rainbow matching in $H$.
This view was applied by Aharoni~\cite{A01} in his proof of Ryser's Conjecture for $r=3$. 
For more applications of this idea, see, e.g., \cite{frankl, AB09, HLS12}.

Conditions of different types are known to guarantee the existence of large rainbow matchings.
For example, in~\cite{penny} a sufficient condition was formulated in terms of domination in an auxiliary graph.
In~\cite{AH00, AB09} conditions were considered in terms of lower bounds on the
size of the largest matching in an auxiliary graph.
For additional problems and results on rainbow matchings, the interested reader is referred to~\cite{WF98, AH00, AB09, KY, AH12}.

In this paper, we study the following setting.
Let $f,t\geq 1$ be integers. A {\em $t$-matching} is a matching of size $t$.
An {\em $(f,t)$-colored} $r$-graph $G=(V,E)$ is an $r$-uniform multi-hypergraph
whose edges are colored in $f$ colors such that every color class is a $t$-matching.
Note that an $r$-subset of $V$ appears in exactly as many of these matchings as its multiplicity in $E$.
We want to determine the smallest number of colors $f$ which guarantees the existence of a rainbow $t$-matching in $G$.
This problem was proposed by Aharoni and Berger~\cite{AB09}, who studied it in the case when 
the $r$-graph is $r$-partite. Formally, let $f(r,t)$ be the largest number $f$ of colors, such
that there exists an $(f,t)$-colored $r$-partite $r$-graph without a rainbow $t$-matching.
If furthermore, each part in this $(f,t)$-colored $r$-partite $r$-graph is required to be of size at most $s$, then
we denote the corresponding extremal value of $f$ by $f_s(r,t)$.
Finally, we denote by $F(r,t)$ the largest value of $f$ such that there exists an $(f,t)$-colored (not necessarily $r$-partite)
$r$-graph with no rainbow $t$-matching.
Obviously,
\begin{equation}
f_s(r,t)\leq f(r,t)\leq F(r,t)
\end{equation}
for every integer $s\geq t$ and $r$.

Aharoni and Berger~\cite{AB09} showed that $f(r,t)\geq 2^{r-1}(t-1)$ for all $r,t > 1$ and proved that 
equality holds for $r=2$ as well as for $t=2$. They also conjectured that their
lower bound is tight in general.
\begin{conjecture}\cite[Conjecture 1.2]{AB09}
\label{conjron}
For every $r,t >1$, $f(r,t)=2^{r-1}(t-1)$.
\end{conjecture}

The lower bound of Aharoni and Berger follows from the following construction.
Let each of the $r$ parts of the vertex set be a copy of $\mathbb{Z}_t$.
For each vector $p \in \{0, 1\}^{r-1}$ define a $t$-matching $M(p)$,
whose $i$-th edge $(1 \leq i \leq t)$ is $\lt(u^{i}_0, u^{i}_1, \ldots, u^{i}_{r-1}\rt)$, 
where $u^{i}_j$ lies in the $j$-th part of the vertex set and
$u_0^i=i$ and $u^{i}_j = i +p(j) \pmod t$ for $j \geq 1$.
To construct a $(2^{r-1}(t-1),t)$-colored $r$-graph without a rainbow $t$-matching,
one can take $t-1$ copies of $M(p)$ for every $p\in \{0,1\}^{r-1}$ as the color classes.

Alon~\cite{noga} realized that the non-existence of rainbow $t$-matchings in the above construction
depends only on the following property of the sequence of the $2^{r-1}$ vectors
$p \in \{0, 1\}^{r-1}$ repeated $t-1$ times: no $t$ of them add up to $0$ in $\mathbb{Z}_t^r$.
Using this he obtained that
\[f(r,t)\geq g(r-1,t)-1,\] where $g(r,t)$ is the well-studied function
denoting the smallest integer $g$ such that any sequence of at least $g$ (not necessarily distinct)
elements of the Abelian group $\mathbb{Z}_t^r$ contains a sub-sequence of exactly $t$
elements whose sum (in  $\mathbb{Z}_t^r$) is zero.
Applying the known lower bounds for the function $g(r-1,t)$, Alon concluded that $f(r,3) > 2.216^r$ for large $r$,
hence Conjecture~\ref{conjron} is false for $t=3$ (as well as for every $t\geq 3$, since $f(r,t) \geq f(r, t-1)$).

Alon also gave a probabilistic construction showing that for large $t$ and all $r$,
\[f(r,t) > 2.71^{r}.\]
As an upper bound he proved that
\begin{equation} \label{eq:noga} f(r,t)\leq F(r,t)\leq \frac{t^{rt}(t-1)}{t!},\end{equation}
which is superexponential in $t$ for fixed $r$.

We substantially improve this estimate, obtaining an upper bound which is polynomial in $t$.
We have two different proof ideas, both giving an estimate of order $t^{2r+1}$. Since 
these ideas may be useful to further improve the bounds and since the proofs are rather short we include  both of them in the paper.
The first one, presented in the next section, gives the following result.
 
\begin{theorem}
\label{theorem}
For arbitrary integers $r,t\geq 2$,
  we have \[f(r,t) \leq F(r,t)<(r+1)^{2r+1}t^{2r+1}.\]
\end{theorem}

Using our second approach, one can slightly reduce this bound on $F(r,t)$ to $\frac{r^r(r+1)^{r+1}}{r!}(t-1)t^{2r}$.
Applying this approach directly to $r$-partite $r$-graphs we are able to improve the leading constant factor
even further.
\begin{theorem}
\label{thm:f<}
For arbitrary integers $r, t\geq 2$,
\[f(r,t)< (r+1)^{r+1}(t-1)t^{2r}.\]
\end{theorem}

In this note we were mainly interested in the case of fixed $r$ and growing $t$.
Comparing the upper bounds on $F(r,t)$, we observe that for very large $r$ the estimate in Theorem~\ref{theorem} is worse than the above mentioned upper bound of Alon~\cite{noga}.
However, using a different argument, one can improve the bound in \eqref{eq:noga}
for essentially all values of $r$.

\begin{theorem}
\label{thm:fixed-t}
For arbitrary integers $r, t\geq 2$,
\[F(r,t)<8^{rt}.\]
\end{theorem}

\noindent
This improves the bound from \eqref{eq:noga}
for all but finitely many pairs $(r,t)$.

\vspace{6mm}
\noindent
{\bf Notation.}
Let $G$ be an $r$-graph. For a subset $\{x_1, \ldots, x_k\}\subseteq V(G) $ we define its
{\em degree} $d(\{x_1, \ldots, x_k\})$, or for short $d(x_1, \ldots, x_k)$,
to be the number of edges containing $\{x_1,\ldots , x_k\}$.
Notice that the degree of an $r$-set is its multiplicity in $E(G)$.
For a subset $S\subseteq V(G)$ of size $k$, we call the {\em link} of $S$ in $G$ the $(r-k)$-graph
whose vertex set is $V(G)\setminus S$ and whose edges are all $\{e\setminus S:\, S\subseteq e\in E(G)\}$
For a vertex $v$, we write  $G(v)$ instead of $G(\{v\})$.

\section{Upper bound for the general case}
In this section, we prove Theorem~\ref{theorem}.
The proof relies on the fact that every multi-hypergraph either has small  covering number,
or contains a large matching. The following technical lemma is a formal consequence of this fact.

\begin{lemma}
\label{lemmacores}
For every $a,b\in \mathbb{N}$ and every (not necessarily uniform) multi-hypergraph $H$ with $\tau(H)\leq a$,
there exists a subset $V'=\{v_1, \ldots, v_c\}\subseteq V(H)$ with  $c=|V'|\leq a$, and a partition
$\lt\{E', E_{v_1},\ldots, E_{v_c}\rt\}$ of the edges of $H$
such that
\begin{itemize}
\item
$|E'|\leq \frac{|E(H)|}{b}$,
\item
$v_i\in e$ for every $i\in [c]$ and every $e\in E_{v_i}$, and
\item
$ |E_{v_i}|>\frac{|E(H)|}{ab}$.
\end{itemize}
\end{lemma}

\proofstart
Let $S$ be a vertex cover of $H$ with $|S|\leq a$. By definition, every edge of $H$ contains some vertex of $S$.
Consider an arbitrary partition $\lt\{E_v:~v\in S\rt\}$
of the edges of $H$ such that for every $v\in S$, the partition satisfies $v\in e$ for every $e\in E_v$.
Denote $E'=\bigcup_{|E_v|\leq \frac{|E(H)|}{ab}} E_v $
and keep those $E_v$ satisfying $|E_v|> \frac{|E(H)|}{ab}$.
Then the second and third properties of the lemma are clearly satisfied.
For the first property, notice that $|E'|\leq a\cdot  \frac{|E(H)|}{ab}=\frac{|E(H)|}{b}$, completing the proof.
\proofend

Armed with this technical result, we are ready to prove Theorem~\ref{theorem}.

\vspace{0.1cm}
\noindent
{\bf Proof of Theorem~\ref{theorem}:} \hspace{2em}
Let $f=(tr+t)^{2r+1}$ and let $G$ be an $(f,t)$-colored $r$-graph.
We aim to prove the existence of a rainbow $t$-matching in $G$.
The main idea of the proof is to show that one of the color classes of $G$ consists of edges
with the property that each of them either has a high multiplicity
or contains a subset whose link has a large vertex cover number.
We then construct our rainbow $t$-matching greedily edge by edge.
We show that both conditions allow us to pick a new edge disjoint from the previously chosen edges
such that the color of the new edge is not present among the colors of the edges chosen before.

As a first step towards this idea, we show that the statement of the theorem holds in case $\tau(G)$ itself is large.
Here and later in the proof, we derive the existence of a large rainbow matching from the fact that 
the vertex cover number is large.

\begin{obs}
\label{largematchinginG}
If $\tau(G) > rt(t-1)$, then $G$ contains a rainbow $t$-matching.
\end{obs}

\noindent
\proofstart
From the assumption it follows by \eqref{eq:tau1} that $\nu(G) \geq \tau(G)/r > t(t-1)$.
Since no color appears more than $t$ times, any maximum matching of $G$ must
contain edges of more than $t-1$ colors, i.e, a rainbow $t$-matching.
\proofend

Hence, we can assume from now on that $\tau(G) \leq rt(t-1)$. We call a non-empty set $S\subseteq V(G)$ a {\em core} of 
$G$ if either $\tau(G(S))> (t-1)(r+1)$ or $|S|=r$ and $d(S)\geq t$.

\begin{lemma}
\label{lemma:disjointcores}
There exists $t$ pairwise disjoint cores.
\end{lemma}

We will show that this follows from the following claim which shows that most of the edges of $G$ contain cores.

\begin{claim}
\label{mostcores}
For every $i\in [ r]$, there exists an $i$-uniform multi-hypergraph $\cF^{(i)}$ on $V(G)$ and a partition 
$\cP^{(i)}=\lt\{E'^{(i)}\rt\}\cup\lt\{\hat{E}^{(i)}\rt\}\cup\lt\{ E_S:~S\in \cF^{(i)}\rt\}$ of the edges of $G$
satisfying
\begin{itemize}
\item[$(i)$]
$\lt|E'^{(i)}\rt|\leq \frac{i|E(G)|}{t(r+1)}$,
\item[$(ii)$]
every $e\in \hat{E}^{(i)}$ contains a core of $G$,
\item[$(iii)$]
for every $S\in \cF^{(i)}$ and every $e\in E_S$, the edge $e$ contains $S$, and
\item[$(iv)$]
every $S\in \cF^{(i)}$ satisfies $|E_S|\geq \frac{|E(G)|}{(t(r+1))^{2i+1}}$.
\end{itemize}
\end{claim}

\proofstart
We prove the claim by induction on $i$.
For $i=1$, we apply Lemma~\ref{lemmacores} with $a=rt(t-1)$ and $b=t(r+1)$ to obtain a partition $\lt\{E', E_{v_1},\ldots, E_{v_c}\rt\}$
of the edges of $G$ satisfying the conditions of the lemma.
We define $E'^{(1)}=E'$ and $\hat{E}^{(1)}$ to be the union of
all the $E_{v_j}$ such that the set $\{v_j\}$ is a core of $G$.
We define the $1$-uniform hypergraph $\cF^{(1)}=\{\{v_j\}:~j\in [c], ~ \{v_j\}\mbox{ is not a core}\}$
and for every $\{v_j\}\in \cF^{(1)}$ we let $ E_{\{v_j\}} = E_{v_j}$.
Then it is easy to check that the corresponding partition $\cP^{(1)}=\lt\{E'^{(1)}\rt\}\cup\lt\{\hat{E}^{(1)}\rt\}\cup\lt\{ E_S:~S\in \cF^{(1)}\rt\}$
satisfies the requirements of the claim.

Let us assume that the assertion of this claim holds for some $i$, $1\leq i < r$, and let $\cF^{(i)}$ and $\cP^{(i)}$ be the corresponding $i$-graph
and the edge partition.
To construct ${\cF}^{(i+1)}$ and $\cP^{(i+1)}$ we initialize $\cF^{(i+1)}=\emptyset$,
$E'^{(i+1)}=E'^{(i)}$ and $\hat{E}^{(i+1)}=\hat{E}^{(i)}$.
Next, for every $S\in \cF^{(i)}$ we distribute the edges in $E_S$ as follows.
Consider the $(r-i)$-graph $G_S^{(i)}:=\lt(V(G)\setminus S, \lt\{e\in\binom{ V(G)}{r-i}:~e\cup S \in E_S\rt\}\rt)$
which is a subgraph of $G(S)$.
Note that $G_S^{(i)}$ is the link of $S$ in the hypergraph $(V(G),E_S)$,
since $S$ is a subset of every $e\in E_S$.
Since $S$ is not a core, the vertex cover number of $G(S)$ and therefore also of $G_S^{(i)}$ is at most $(t-1)(r+1)$.
Applying Lemma~\ref{lemmacores} with $a=(t-1)(r+1)$ and $b=t(r+1)$ to $G_S^{(i)}$
we obtain a partition $ \lt\{E', E_{v_1}, \ldots, E_{v_c}\rt\}$ of $E\lt(G_S^{(i)}\rt)$.
This partition satisfies
\begin{itemize}
\item
$|E'|\leq  \frac{\lt|E\lt(G_S^{(i)}\rt)\rt|}{b} = \frac{|E_S|}{t(r+1)}$
\item
$v_j\in e$ for every $j\leq c$ and every $e\in E_{v_j}$, and
\item
$\lt|E_{v_j}\rt|\geq  \frac{\lt|E\lt(G_S^{(i)}\rt)\rt|}{ab}  = \frac{|E_S|}{t(t-1)(r+1)^2}> \frac{|E(G)|}{(t(r+1))^{2i+3}}$ for every $j\leq c$.
\end{itemize}
For the last inequality we used that by induction $|E_S|\geq \frac{|E(G)|}{(t(r+1))^{2i+1}}$.

For every $e\in E'$ we add $e\cup S$ to $E'^{(i+1)}$.
For the rest,
if $S\cup \{v_j\}$ is a core, then for every edge $e\in E_{v_j}$ we add $e\cup S$ from $E_S$ to $\hat{E}^{(i+1)}$.
If $S\cup \{v_j\}$ is not a core, we add $S\cup \{v_j\}$ to $\cF^{(i+1)}$  and for every edge $e\in E_{v_j}$ we add $e\cup S$ from $E_S$ to
$E_{S\cup \{v_j\}}$.
Note that this way, several copies of $S\cup \{v_j\}$ could appear in $\cF^{(i+1)}$ simply because there might be
several ways to split an $i$-set into a $1$-set and an $(i-1)$-set.
This is why we allow $\cF^{(i)}$ to be a multiset. Notice however that each copy of $S\cup \{v_j\}$ in $\cF^{(i+1)}$ is assigned to a distinct
edge multiset $E_{S\cup \{v_j\}}$ (which is disjoint from the others).

After distributing all the edges of each $E_S$ we obtain a partition $\cP^{(i+1)}$ of the edges of $G$
with a corresponding $(i+1)$-graph $\cF^{(i+1)}$.
By construction every edge from $\hat{E}^{(i+1)}$ contains a core.
We also made sure that for every $Q\in \cF^{(i+1)}$,
every $e\in E_Q$ contains $Q$. Furthermore, $|E_Q|\geq \frac{|E(G)|}{(t(r+1))^{2i+3}}$, since 
$E_Q$ consists of sets that are the union of some $S\in {\cF^{(i)}}$ and the members of some $E_{v_j}$, and therefore $|E_Q|=|E_{v_j}|$.
Finally, note that $E'^{(i+1)}$ contains edges from $E'^{(i)}$ as well as at most a $\frac{1}{t(r+1)}$-fraction of each $E_S$. Therefore
we have 
\[|E'^{(i+1)}|\leq \frac{i|E(G)|}{t(r+1)}+ \frac{\sum |E_S|}{t(r+1)}\leq 
\frac{i|E(G)|}{t(r+1)}+\frac{|E(G)|}{t(r+1)}=
\frac{(i+1)|E(G)|}{t(r+1)}.\]
\proofend

We are now ready to prove Lemma~\ref{lemma:disjointcores} using Claim~\ref{mostcores}.

\noindent
{\bf Proof of Lemma~\ref{lemma:disjointcores}:} \hspace{2em}
Applying Claim~\ref{mostcores} with $i=r$, we have that the only edges not containing a core of $G$ are all in $E'^{(r)}$.
Indeed, the edges from $\hat{E}^{(r)}$ contain a core of $G$ by the part $(ii)$ of the claim.
Also every $S\in \cF^{(r)}$ is a core, since $|S|=r$ and by part $(iv)$ the multiplicity of $S$ in $E(G)$ is at least 
\[|E_S|>\frac{ |E(G)|}{(t(r+1))^{2r+1}}=
\frac{f\cdot t}{(t(r+1))^{2r+1}} \geq t.\]
On the other hand, by part $(i)$ of Claim ~\ref{mostcores} $\lt|E'^{(r)}\rt|< |E(G)|/t=f$. Thus there exists a color class containing no edges from $E'^{(r)}$.
Every edge of this color class contain a core which gives $t$ disjoint cores.
\proofend

Let $S_1,\ldots, S_t$ be $t$ disjoint cores from Lemma~\ref{lemma:disjointcores}.
To finish the proof of Theorem~\ref{theorem} we iteratively find a rainbow matching $\{e_1, \ldots , e_t\}$, such that
$e_j\supseteq S_j$ as follows. Assume that for some $j=0, 1, \ldots , t-1$ we already have a rainbow $j$-matching $\{e_1, \ldots,e_{j}\}$ whose edges are disjoint from the $(t-j)$ 
cores $S_{j+1},\ldots, S_t$. We show how to find an edge $e_{j+1}\in E(G)$ such that $e_1, \ldots, e_{j+1}$ is a rainbow $(j+1)$-matching,
and $e_{j+1}$ is disjoint from every core $S_{j+2},\ldots, S_t$. Hence, in the end, we have a rainbow $t$-matching.

First consider the case when $S_{j+1}$ has size less than $r$. Then we have
$\tau(G (S_{j+1})) > (t-1)(r+1)$. Let $U=\bigcup_{i=1}^je_i \cup \bigcup_{\ell\geq j+2}S_\ell$. Then $|U| \leq (t-1)r$.
Note that for any hypergraph $H$ and a subset of vertices $W\subseteq V(H)$ 
the number of edges of $H$ disjoint from $W$ is at least $\tau (H)-|W|$.
Thus, taking $H=G(S_{j+1})$ and $W=U$, we have that the number of edges in
$G(S_{j+1})$ disjoint from $U$ is at least
$\tau(G(S_{j+1})) - |U| >  t-1$.
Therefore, by definition of $G(S_{j+1})$ there exist $t$ edges $g_1, \ldots ,g_t$ in $G$ which contain 
$S_{j+1}$  and are disjoint from the
edges of $e_1, \ldots, e_j$ and from the cores $S_{j+2},\ldots, S_t$.
Since the edges $g_1, \ldots ,g_t$ are all pairwise intersecting, they have $t$ distinct colors.
Thus one of the $g_\ell$ has a color that is different from all the colors of the edges $e_1, \ldots, e_j$.
This edge, which we denote by $e_{j+1}$, satisfies the requirements of the iteration.

In the second case when $|S_{j+1}|=r$,  by the definition of a core we have that $S_{j+1}$ is an edge of $G$ with multiplicity at least $t$.
Then there is a color of this edge which is distinct from colors of $e_1, \ldots, e_j$. Choosing
$e_{j+1}$ to be $S_{j+1}$ with this color satisfies the requirements of the iteration and completes the proof of the theorem.
\proofend

\section{Upper bound for the $r$-partite case}
\label{rest}

In this section we slightly improve the bound on $f(r,t)$ from Theorem~\ref{theorem} by
giving a completely different proof. This proof can also be easily adapted to not necessarily $r$-partite $r$-graphs
(see below).

When thinking about the smallest value $f$ forcing an $r$-graph $G$ to contain a rainbow $t$-matching,
one is necessarily confronted with the question about the structure of the extremal examples.
It is somewhat intuitive to expect that the extremal $r$-graphs are dense in the sense that they should not contain unnecessary vertices.
The following lemma shows that this intuition is indeed correct.
It states that for sufficiently large values of $t$ there exist nearly optimal constructions with all $r$ parts being of size not much
larger than $rt^2$.

Recall that by $f_s(r,t)$ we denote the maximum integer $f$ for which there exists an $(f,t)$-colored $r$-partite $r$-graph
with parts of size at most $s$ not containing a rainbow $t$-matching.

\begin{lemma}
\label{dense}
For every $r$, $t\geq 1$ and $s\geq t$,
$f_s(r,t)> \lt(1-\frac{rt^2}{s}\rt) f(r,t)$.
\end{lemma}

\proofstart
Let $G$ be an $(f,t)$-colored $r$-partite $r$-graph  with $f=f(r,t)$,
which contains no rainbow $t$-matching.
Starting from $G$, we iteratively construct an $(f_0,t)$-colored
$r$-partite $r$-graph $G_0$ with parts of size at most $s$ and with $f_0> \lt(1-\frac{rt^2}{s}\rt)f$,
such that $G_0$ also contains no rainbow $t$-matching.
The main idea of the proof is that contracting two vertices from the same part of $V(G)$
does not create a rainbow $t$-matching. In order to keep most of the colors to be $t$-matchings,
one just have to make sure that the two contracted vertices do not appear simultaneously on two 
edges of many color classes.

Given any $r$-partite $r$-graph $H$ with parts $V_1 \cup \cdots \cup V_r$, whose edges are colored such that all color classes are $t$-matchings,
consider the following auxiliary $t$-graph $H'$. It has the same vertex set $V(H')=V(H) = V_1 \cup \cdots \cup V_r$
and we put a $t$-edge $\{ v_1, \ldots, v_t\} \subseteq V_i$ into $H'$ if and only if
the edges of some color class of $H$ intersect $V_i$ exactly in  $\{ v_1, \ldots, v_t\}$.
Notice that $H'$ is the union of $r$ vertex-disjoint $t$-graphs whose vertex sets are the parts of $V(H)$.
Furthermore, note that for every vertex $v\in V(H)$, we have $d_H(v)=d_{H'}(v)$.
To make the notation consistent, the prime-sign always denotes the auxiliary $t$-graph.

Starting with $G$, we iteratively perform the following transformation of our $r$-graph.
Suppose that we are currently dealing with an $(\hat{f}, t)$-colored $r$-partite $r$-graph $\hat{G}$
not containing a rainbow $t$-matching.
Choose arbitrarily a part $V_k$ with $|V_k|>s$, if such part exists.
Take two distinct vertices $x$ and $y$ in $V_k$
whose degree in the corresponding auxiliary $t$-graph $\hat{G}'$ is smallest among all the pairs of distinct vertices from $V_k$.
By double counting the sum of the degrees of all $2$-subsets of $V_k$,
we obtain
\begin{equation*}
d_{\hat{G}'}(x,y) \binom{|V_k|}{2}\leq\sum_{u,w \in V_k,\, u\neq w}d_{\hat{G}'}(u,w)= \hat{f}\binom {t}{2},
\end{equation*}
since each color class contributes ${t\choose 2}$ to the sum.
Hence, using the fact that $|V_k|\geq s\geq t$, we obtain
\begin{equation}
\label{s<}
d_{\hat{G}'}(x,y) <\frac{\hat{f}t^2 }{|V_k|^2}.
\end{equation}

Delete from $\hat{G}$ every edge of a color $c$ for which there exist edges $e_1, e_2$ both of color $c$
such that $x\in e_1$ and $y\in e_2$.
Notice that by~\eqref{s<} the total number of such colors is at most $\frac{\hat{f}t^2}{|V_k|^2}$.
In the resulting $r$-graph, there exists no color class containing $x$ and $y$ simultaneously.
Replace every appearance of $x$ and $y$ in every remaining edge by a new auxiliary vertex $z\not \in V(\hat{G})$,
thus reducing the size of $V_k$ by one. Notice that the new $r$-graph has no rainbow $t$-matching
as well, and every color class in its edge-coloring is still a $t$-matching.

Iterating this transformation until all parts have size at most $s$,
we obtain an $(f_0,t)$-colored $r$-partite $r$-graph $G_0$ with parts of size at most $s$
not containing a rainbow $t$-matching .
For the number of colors $f_0$ in the coloring of $E(G_0)$, 
observe by \eqref{s<} that during the iterations we deleted at most
\[r\sum_{\ell> s}\frac{ft^2 }{\ell^2}<frt^2\int_{\ell=s}^{\infty}\frac{1}{\ell^2}d\ell=\frac{frt^2}{s}\]
color classes from $G$.
Hence $E(G_0)$ contains more than $\lt(1-\frac{rt^2}{s}\rt)f$ color classes, completing the proof of the lemma.
\proofend

With Lemma~\ref{dense} in our hand, it suffices to give an upper bound on $f_s(r,t)$ to obtain an upper bound on $f(r,t)$.

\begin{lemma}
\label{f^-}
For every $r,t\geq 1$ and $s\geq t$, $f_s(r,t)\leq (t-1)s^{r}$.
\end{lemma}

\proofstart
Let $f:=(t-1)s^r+1$ and let $G$ be an $(f,t)$-colored $r$-partite $r$-graph with parts of size at most $s$.
We prove that $G$ contains a rainbow $t$-matching.
Call an $r$-set containing one vertex from every part of $V(G)$ {\em bad} if its multiplicity in $E(G)$ is less than $t$.
Since there are only at most $s^r$ bad $r$-sets (as there are no more eligible $r$-sets),
and each of them appears in at most $t-1$ colors,
there must be a color class in $G$ which does not contain any bad set.
This color class is a $t$-matching $M$ and each of its edges has multiplicity at least $t$.
Thus we can construct greedily a rainbow $t$-matching by picking edges of $M$ colored with distinct colors. 
\proofend

Theorem~\ref{thm:f<} now follows immediately from
the above two lemmas.

\noindent
{\bf Proof of Theorem~\ref{thm:f<}:} \hspace{2em}
By  Lemmas~\ref{dense} and~\ref{f^-} we have for every $r,t\geq 1$ and $s\geq t$ that
\[ f(r,t) < \frac{1}{1-\frac{rt^2}{s} } \cdot (t-1)s^r= \frac{(t-1)s^{r+1}}{s-rt^2}.\]
Substituting $s=t^2(r+1)$ completes the proof.
\proofend

One can also use these ideas in the non-partite setting. Indeed, given any $(f,t)$-colored
$r$-graph one can find a pair of vertices whose contraction destroys only at most 
$\frac{r^2t^2}{|V(G)|^2}f$ color classes. Therefore as we showed above one can reduce the vertex set of $G$ to be of size 
$s$ and still have more than $(1-\frac{r^2t^2}{s})f$ $t$-matchings. Once the ground set has size $s$, there are at most 
$(t-1){s \choose r}\leq (t-1)s^r/r!$ color classes in which some edge has multiplicity at most $t-1$. Therefore if the the number of 
colors is greater than that, there is a color class with all the edges having multiplicity at least $ t$ and hence there is a 
rainbow $t$-matching. Combining these arguments and choosing $s=r(r+1)t^2$ shows that
$F(r,t)< \frac{r^r(r+1)^{r+1}}{r!}(t-1)t^{2r}$.
The rest of the details are very similar to the proof in the $r$-partite case 
and are omitted.

\section{The case of small $t$ and large $r$}

The following theorem gives the recursion which implies Theorem~\ref{thm:fixed-t}.
Its proof combines the approach from \cite{noga} with some additional ideas.

\begin{theorem}
\label{improve}
For every $r,t\geq 2$,
\begin{equation}
\label{smallt-larger}
F(r,t) \leq \frac{2^{rt}}{\binom{t}{\lt\lceil\frac{t}{2}\rt\rceil}}  \lt(F\lt(r,\lt\lceil\frac{t}{2}\rt\rceil\rt)+\lt\lfloor\frac{t}{2}\rt\rfloor\rt).
\end{equation}
In particular, $F(r,t)<8^{rt}$.
\end{theorem}

\proofstart
Let $G$ be an $(f,t)$-colored $r$-graph with 
\[f> \frac{2^{rt}}{\binom{t}{\lt\lceil\frac{t}{2}\rt\rceil}}  \lt(F\lt(r,\lt\lceil\frac{t}{2}\rt\rceil\rt)+\lt\lfloor\frac{t}{2}\rt\rfloor\rt).\]
Color each {\em vertex} of $G$ independently, uniformly at random with black and white.
We say that a color class 
{\em survives} this procedure if all the vertices in exactly $\lt\lfloor\frac{t}{2}\rt\rfloor$ of its $t$ edges become black,
and all the vertices in its remaining $\lt\lceil\frac{t}{2}\rt\rceil$ edges become white.
Note that each color class survives with probability
$\binom{t}{\lt\lceil\frac{t}{2}\rt\rceil}2^{-rt}$, because specifying the $\lt\lceil\frac{t}{2}\rt\rceil$ of its $t$ edges whose vertices should
become white fully determines the color of all the $rt$ vertices of this $t$-matching.
Thus, by linearity of expectation,
there exists a coloring of the vertices of $G$ such that
\[f'\geq \frac{\binom{t}{\lt\lceil\frac{t}{2}\rt\rceil}}{2^{rt}} \cdot f
>
F\lt(r,\lt\lceil\frac{t}{2}\rt\rceil\rt)+\lt\lfloor\frac{t}{2}\rt\rfloor
\]
of the $f$ color classes survive.

Consider the $r$-graph remaining after deleting edges
of all the color classes that do not survive.
Then the above argument guarantees the existence of an $r$-graph $G'$ with the following properties:
\begin{itemize}
\item
the vertex set of $G'$ is colored white and black;
\item
the edges of $G'$ are colored with $f'$ colors such that every color class is a $t$-matching;
\item
every color class consists of $\lt\lfloor\frac{t}{2}\rt\rfloor$ edges containing only black vertices
and $\lt\lceil\frac{t}{2}\rt\rceil$ edges containing only white vertices.
\end{itemize}

Let $G_b$ be the $r$-graph obtained from $G'$ by deleting all white vertices (and all the edges that contain them).
By construction, $G_b$ is an $(f',\lt\lfloor\frac{t}{2}\rt\rfloor)$-colored $r$-graph.
Since $f'> F\lt(r,\lt\lceil\frac{t}{2}\rt\rceil\rt)\geq F\lt(r,\lt\lfloor\frac{t}{2}\rt\rfloor\rt)$, we know that $G_b$ contains a rainbow 
$\lt\lfloor\frac{t}{2}\rt\rfloor$-matching $M_b$.

Now let $G_w$ be the $r$-graph obtained from $G'$ by deleting all black vertices (and the edges that contain them)
as well as all edges contained in the same color class with one of the edges from $M_b$.
Then $G_b$ is an $(f'-\lt\lfloor\frac{t}{2}\rt\rfloor, \lt\lceil\frac{t}{2}\rt\rceil)$-colored
$r$-graph. Since $f'-\lt\lfloor\frac{t}{2}\rt\rfloor >  F\lt(r,\lt\lceil\frac{t}{2}\rt\rceil\rt)$,
we know that $G_w$ contains a rainbow $\lt\lceil\frac{t}{2}\rt\rceil$-matching $M_w$.
Now the union of the two matchings $M_b$ and $M_w$ is a rainbow $t$-matching in $G$,
proving that $f>F(r,t)$.

Finally, using (\ref{smallt-larger}), together with the  obvious fact that
\[F(r,1)=0<8^{r\cdot 1},\]
one can check by induction that $F(r,t)<8^{rt}$.
\proofend

\section{Concluding remarks and open problems}

Although the conjecture of Aharoni and Berger was refuted by Alon,
we still believe that for fixed $r$, the function $f(r,t)$ grows linearly in $t$.

\begin{conjecture}
\label{conjecture}
For every $r$ there exists a constant $c_r$ such that $f(r,t) \leq c_rt$ for all $t$.
\end{conjecture}

\noindent
Since currently we can only prove polynomial bound whose exponent depends on $r$, it would be even interesting to prove that there exists a function $c_r$ depending only on $r$ 
and an absolute constant $b$, such that $f(r,t) \leq c_r t^b$.

As we already discussed in Section~\ref{rest}, we do believe that extremal configurations for this problem should not have too many vertices.
In particular, it would be interesting to decide whether there are $(f,t)$-colored $r$-graphs on $O(rt)$ vertices with 
$f=\Theta\big(F(r,t)\big)$ and no rainbow $t$-matching. Similarly, one can ask whether
$f(r,t)=\Theta\big(f_s(r,t)\big)$ for some $s=O(t)$.

Another natural question concerns the value of $f(r,t)$ when $t$ is fixed and $r$ grows. We know that
it grows exponentially in $r$ and for large $t$ we have
\[2.71^r < f(r,t) < 8^{tr}.\]
It would be interesting to determine whether $f(r,t)$ can be upper bounded by $\alpha_t {\beta}^r$ for some absolute constant $\beta$ and some 
function $\alpha_t$ depending on $t$.

\vspace{0.2cm}

\noindent
{\bf Acknowledgment.} We would like to thank Yury Person for fruitful discussions at the early stages of this project.

\end{document}